\documentclass[11pt]{article}
\usepackage{amssymb}
\usepackage{amsmath}
 \usepackage{titlesec}
\usepackage{fancyhdr}
\usepackage{geometry}
\usepackage{mathrsfs}
\usepackage[english]{babel}
\usepackage{amsthm}
\usepackage{mathtools}

\newtheorem{mainthm}{Theorem}

\newtheorem{theorem}{Theorem}
\usepackage{lipsum}
\newtheorem{lemma}{Lemma}
\newtheorem{corollary}{Corollary}
\usepackage[hang,flushmargin]{footmisc} 
\numberwithin{equation}{section}
\usepackage[skip=10pt plus1pt, indent=25pt]{parskip}
\newcommand\blfootnote[1]{%
  \begingroup
  \renewcommand\thefootnote{}\footnote{#1}%
  \addtocounter{footnote}{-1}%
  \endgroup
}

\usepackage[dvipsnames]{xcolor}
\usepackage{ifpdf}
\ifpdf
\usepackage[pdftex]{graphicx}
\usepackage{epstopdf}
\usepackage{tikz}
\usepackage{hyperref}
\usepackage{doi}
\usepackage[nameinlink]{cleveref}

\usepackage[utf8]{inputenc}
\else
  \usepackage[dvips]{graphicx}
\fi
\newcommand\myshade{85}
\colorlet{mylinkcolor}{blue}
\colorlet{mycitecolor}{blue}
\colorlet{myurlcolor}{Aquamarine}
\hypersetup{linkcolor  = mylinkcolor,
	citecolor  = mycitecolor,
	urlcolor   = myurlcolor!\myshade!black,
	colorlinks = true,}

\title{ \Large\bfseries{A bound for Stieltjes constants of the Riemann zeta function}}
\author{{\small Omprakash Atale${}^{1}$,\,\,\,Suraj Singh Khurana${}^{2}$ }\\[1pt]\\{\small${}^{1,2}$ Department of Mathematics,} \\{\small SRM University AP, Amaravati-522240, }\\{\small  Andhra Pradesh, India.}}
\date{}

\begin{document}

\maketitle
\vspace{1em}
\begin{abstract}
The coefficients appearing in the Laurent expansion of the Riemann zeta function are known as the Stieltjes constants and are denoted by $\gamma_n$. Several authors, including Berndt, Matsuoka, Williams and Zhang, and Saad Eddin, have investigated upper bounds for $|\gamma_n|$, with the best known result due to Pauli and Saidak \cite{MR4671618}. In 2003, based on numerical evidence, Kreminski \cite{MR1972742} conjectured that Matsuoka's bound \cite{MR827790} could be considerably improved. Motivated by this conjecture, we employ the saddle-point method to derive an asymptotic formula for the Stieltjes constants, with the saddle point expressed explicitly in terms of the Lambert $W$-function. The resulting formula yields an asymptotic upper bound of the correct order of magnitude. Moreover, the resulting bound is consistent with the conjecture of Farr, Pauli, and Saidak \cite{MR4280128}.\blfootnote{${}^{1}$E-mail:- atale.om@gmail.com\\ ${}^{2}$E-mail:- suraj.singh.khurana@gmail.com }
\end{abstract}
\vspace{1em}

\section{Introduction}
The Riemann zeta function is defined by the Dirichlet series 
\begin{equation}
    \zeta(s):=\sum_{n=1}^{\infty}\frac{1}{n^s},
\end{equation}
where $s=\sigma+it$, and the series converges for $\sigma>1$ \cite[\S 1.1]{MR882550}. The function $\zeta(s)$ extends to a meromorphic function on $\mathbb{C}$ with a simple pole at $s=1$, and satisfies the functional equation
\begin{equation}
    \zeta(s)=2^s\pi^{s-1}\sin\bigg(\frac{1}{2}\pi s\bigg) \Gamma(1-s)\zeta(1-s).
\end{equation}
The Laurent series expansion of $\zeta(s)$ around $s=1$ is given by \cite[p. 325]{MR69209}
\begin{equation}
    \zeta(s)=\frac{1}{s-1}+\sum_{n=0}^\infty\frac{(-1)^n}{n!}\gamma_n(s-1)^n,
\end{equation}
where
\begin{equation}
\gamma_n=\lim_{N\to\infty}\bigg(\sum_{k=1}^N\frac{\log^n k}{k}-\frac{\log^{n+1}N}{n+1}\bigg).
\end{equation}
The coefficients $\gamma_n$ are known as the Stieltjes constants and $\gamma_0$ is the well-known Euler's constant. As recently shown in \cite{MR3472126, MR3853435, MR4120648, MR2728984, MR4383337}, these constants have a number of interesting and unexpected applications in the theory of the zeta function. These constants have been studied by many authors; see, for example, \cite{MR584816, MR662837, MR657478, MR1302533, MR1277891, MR161842, MR2280843}. 

In 1955,  Briggs \cite{MR76858} established the first explicit upper bound for the magnitude $|\gamma_n|$. This initial result was subsequently improved by Berndt in 1972 \cite{MR289431}, Israilov in 1981 \cite{MR662837} and Williams and Zhang \cite{MR1302533} in 1994. In 1985, Matsuoka \cite[Theorem 6]{MR827790} proved that    
\begin{equation}
    |\gamma_n|<10^{-4}(\log n)^n
\end{equation}
for $n>9$, and in 1989, he obtained an asymptotic expansion for $\gamma_n$ \cite{MR1001731}. In 2003, on the basis of numerical evidence, Kreminski \cite[\S 1]{MR1972742} conjectured that the above bound could be strengthened considerably. Matsuoka's bound remained the sharpest known upper bound for $|\gamma_n|$ for more than two decades. Several attempts have been made to improve Matsuoka's bound, but only a few have succeeded. The next improvement came in 2011, when Adell \cite{MR2774545}, using a probabilistic approach, based on a differential calculus for the gamma process, obtained an explicit upper bound for $|\gamma_n|$, valid for $n>3$, whose order of magnitude is
\begin{equation}
    \exp\bigg(n\log\log n-\frac{n\log\log n}{\log n}(1+O(1))\bigg).
\end{equation}
In the same year, Knessl and Coffey \cite{MR2728984} obtained an asymptotic formula that approximates the Stieltjes constants remarkably well for large $n$. 

The most recent refinement is due to Pauli and Saidak in 2024 \cite{MR4671618}. They proved that 
\begin{equation}
    |\gamma_n|<n^2+\frac{3}{4}n^2\log n
\big|e^{n\log w(n)-\frac{n}{w(n)}}\big|,\label{PS}
\end{equation}
where $n>2\pi$, $w(n):=W_0\big(i\frac{n}{2\pi}\big)$ and $W_0$ is the principal branch of the Lambert $W$ function (cf. \eqref{WAE}). Farr, Pauli and Saidak \cite[p. 21]{MR4280128} conjectured that
\begin{equation}
    |\gamma_n|<2\big|e^{n\log w(n)-\frac{n}{w(n)}}\big|,\label{CON}
\end{equation}
so \eqref{PS} differs from the conjectured bound only by a factor of order $n^2\log n$.  From a computational point of view, the formulas discussed above still have some limitations \cite[\S 1.1]{MR3925487}. In fact, none of the known upper bounds for $|\gamma_n|$ discussed above gives the correct order of magnitude. Motivated by this, we prove Theorem 1 for $\gamma_n$.

In what follows, let
\begin{equation}
    T(n):=\log\bigg(\bigg(n+\frac{3}{2}\bigg)\frac{1}{2\pi}\bigg)-i\frac{\pi}{2}.
\end{equation}

\begin{theorem}\label{Thm1}
    For all sufficiently large integers $n$,
    \begin{equation}\label{MT}
        \gamma_n= \sqrt{\frac{2}{\pi}}e^{R(n)}\bigg(\cos J(n)+O\bigg(\frac{\log^{2/3 }n}{n^{1/3}}\bigg)\bigg)
    \end{equation}
where $R(n)$ and $J(n)$, defined by \eqref{R} and \eqref{J}, satisfy
\begin{equation}
    \begin{aligned}
        R(n)=\bigg(n+&\frac{3}{2}\bigg)\bigg(\Re\log\lambda(n)-\Re\frac{1}{\lambda(n)}\bigg)-\frac{1}{2}\log\bigg|\lambda^2(n)+\lambda(n)+\frac{2\lambda^3(n)}{(2\sqrt{3}n+3\sqrt{3})^2}\bigg|\\&-\frac{1}{2}\log\bigg(n+\frac{3}{2}\bigg)+\frac{11+4\Re\,\lambda(n)}{24n+36}+\frac{1}{(2\sqrt{2}n+3\sqrt{2})^2}+\log 2\pi+O\bigg(\frac{\log^3 n}{n^3}\bigg),
    \end{aligned}
\end{equation}
\begin{equation}
    \begin{aligned}
        J(n)=\bigg(n+\frac{3}{2}\bigg)&\bigg(\Im\log \lambda(n)-\Im\bigg(\frac{1}{\lambda(n)}\bigg)\bigg)+\frac{\Im\lambda(n)}{6n+9}\\&-\frac{1}{2}\arg\bigg(\lambda^2(n)+\lambda(n)+\frac{2\lambda^3(n)}{(2\sqrt{3}n+3\sqrt{3})^2}\bigg)+O\bigg(\frac{\log^3n}{n^3}\bigg),
    \end{aligned}
\end{equation}
respectively, and $\lambda(n)$, defined by \eqref{lambda}, satisfies
 \begin{equation}
        \lambda(n)=W_0(e^{T(n)})-\frac{W_0(e^{T(n)})^3}{\big(2\sqrt{3}n+3\sqrt{3}\big)^2\big(W_0(e^{T(n)})+1\big)}+O\bigg(\frac{\log^4 n}{n^4}\bigg).
    \end{equation}
\end{theorem}
The following corollary is an immediate consequence of \autoref{Thm1} and the triangle inequality.
\begin{corollary}
    For all sufficiently large integers $n$,
    \begin{equation}\label{NB}
        |\gamma_n|\leq \sqrt{\frac{2}{\pi}}e^{R(n)}\bigg(1+O\bigg(\frac{\log^{2/3 }n}{n^{1/3}}\bigg)\bigg).
    \end{equation}
\end{corollary}

A numerical comparison of the main terms of \eqref{MT} and \eqref{NB} with the bounds \eqref{PS} and \eqref{CON} is given in \hyperref[fig1]{Figure \ref{fig1}} and \hyperref[tab:comparison]{Table \ref{tab:comparison}}.

\subsection{Asymptotic expansion of $\gamma_n$}
Our main results rest on an asymptotic expansion of $\gamma_n$ obtained by Matsuoka, and we begin by recalling it. Matsuoka writes $\gamma_n$ as a contour integral, using Cauchy's integral formula and the functional equation of $\zeta(s)$, and then evaluates this integral by the saddle-point method (cf. \cite{MR2139829, MR4671618, MR3925487}). We now give a brief outline of Matsuoka's method. In what follows, let
\begin{equation}
    \phi(z):=-(n+1)\log z-z\log(2\pi i)+\log \Gamma(z)
\end{equation}
and
\begin{equation}
    \Omega_n:=\big\{z=x+iy\in\mathbb{C}: \sqrt{n}<x<n,\,0<y<x\big\}.
\end{equation}
The relevant saddle point lies in this region. Matsuoka shows that, for all sufficiently large integers $n$, the function $\phi(z)$ has the following properties in $\Omega_n$ \cite[Lemmas 1--4]{MR1001731}:
\begin{enumerate}
    \item The equation $z\phi'(z)=0$  has a unique solution in the domain $\Omega_n$, which we denote by $z_0=a+ib$. We will refer to this equation as the saddle-point equation.
    \item The real and imaginary parts $a=a(n)$ and $b=b(n)$ depend only on $n$ and satisfy
    \begin{equation}\label{a(n)}
        a(n)=\frac{n}{\log n}+O\bigg(\frac{n\log\log n}{\log^2 n}\bigg),
    \end{equation}
    \begin{equation}
        b(n)=\frac{\pi}{2}\frac{n}{\log^2 n}+ O\bigg(\frac{n\log\log n}{\log^3 n}\bigg).
    \end{equation}
    \item The function $g(y):=\Re\,\phi(a+iy)$ is strictly increasing on $0\leq y\leq b$, attains its maximum at $y=b$, and is strictly decreasing on $b\leq y\leq a$.
    \item If $\Delta:=\sqrt{a}\log a$, then
    \begin{equation}
        g(b)-g(b\pm\Delta)>\frac{1}{3}\log^ 3a.
    \end{equation}
    \item The functions $g(y)$ and $f(y):=\Im\,\phi(a+iy)$ satisfy \footnote{For the real functions $g$ and $f$, primes denote derivatives with respect to $y$; for $\phi$, they denote complex derivatives.}
    \begin{equation}\label{g''}
        g''(b)=-\frac{\log a}{a}+(\log 2\pi-1)\frac{1}{a}+O\bigg(\frac{1}{a\log a}\bigg),\quad\mathrm{and}
    \end{equation}
    \begin{equation}
        f''(b)=\frac{\pi}{a}+O\bigg(\frac{1}{a\log a}\bigg).
    \end{equation}
\end{enumerate}
It is clear from the above properties that $z_0$ is a saddle point of $\phi(z)$. From the definition of $\phi(z)$, one has
\begin{equation}
    g(y)=-\frac{1}{2}(n+1)\log\big(a^2+y^2\big)-a\log 2\pi
+\frac{1}{2}\pi y+\Re\log\Gamma(a+iy),\quad\mathrm{and}
\end{equation}
\begin{equation}
    f(y)=-(n+1)\arctan\bigg(\frac{y}{a}\bigg)-y\log 2\pi-\frac{1}{2}\pi a+\Im\log\Gamma(a+iy).
\end{equation}
The expressions for $\Re\log\Gamma(a+iy)$ and $\Im\log\Gamma(a+iy)$ appearing above can be evaluated further using Stirling's formula \cite[Eqn. (1.4.25)]{MR174795}. 

Let $\left\{h_j\right\}_{j\geq 0}$ be a sequence of numbers characterized by the expansion
\begin{equation}
    \sum_{j=0}^\infty h_j(y-b)^j=\exp\bigg(\phi(a+iy)-\phi(a+ib)+\frac{1}{2}\phi''(a+ib)(y-b)^2\bigg).
\end{equation}
Moreover, define the sequences $\left\{u_j\right\}_{j\geq 0}$ and $\left\{v_j\right\}_{j\geq 0}$ by $u_j:=\Re\,h_j$ and $v_j:=\Im\,h_j$ respectively. Matsuoka's main result is the following.

\begin{mainthm}[Matsuoka \cite{MR1001731}]\label{Matsuoka}
    Let $N$ be a nonnegative integer. If $n>c\,\Gamma\big(N/3+17/6\big)$, where $c$ is a sufficiently large constant, then
    \begin{equation}
        \begin{aligned}
            \gamma_n=\frac{1}{\pi}&n!e^{g(b)}\sum_{j=0}^{N}|h_{2j}|2^{j+1/2}\Gamma\bigg(j+\frac{1}{2}\bigg)\big(g''(b)^2+f''(b)^2\big)^{-j/2-1/4}\\&\times\cos\bigg\{f(b)-\bigg(j+\frac{1}{2}\bigg)\arctan\bigg(\frac{f''(b)}{g''(b)}\bigg)+\arctan\bigg(\frac{v_{2j}}{u_{2j}}\bigg)\bigg\}\\&+O\bigg(\Gamma\bigg(N+\frac{3}{2}\bigg)2^Nn!e^{g(b)}n^{-N/3+1/6}\log^{2N/3-1/3}n\bigg).
        \end{aligned}
    \end{equation}
\end{mainthm}
We denote by $G(n)$ and $F(n)$ the quantities $g(b)$ and $f(b)$, respectively, thereby emphasizing the fact that the quantities depend only on $n$. Thus,
\begin{equation*}
    G(n):=\Re\,\phi(a(n)+ib(n))\quad\mathrm{and}\quad F(n):=\Im\,\phi(a(n)+ib(n)).
\end{equation*}
The simplest case of \autoref{Matsuoka} is $N=0$, which gives
\begin{equation}
    \gamma_n=\bigg(\frac{2}{\pi}\bigg)^{1/2}n!\frac{n^{1/2}}{\log n}e^{G(n)}\bigg(\cos F(n)+O\bigg(\frac{\log\log n}{\log n}\bigg)\bigg)
\end{equation}
with
\begin{equation}
    F(n)=-\frac{\pi}{2}\frac{n}{\log n}+O\bigg(\frac{n\log\log n}{\log^ 2n}\bigg).
\end{equation}
Thus the size of $|\gamma_n|$ is controlled by $G(n)$, while the sign of $\gamma_n$ is controlled by $F(n)$. However, $G(n)$ and $F(n)$ are defined only implicitly, through the saddle point $z_0$, and the absolute errors in the above estimates for $a(n)$ and $b(n)$ grow without bound as $n\to\infty$. 

\section{Main results} 
The Lambert $W$ function $W(z)$ is defined by the equation $W(z)e^{W(z)}=z$. The function $W(z)$ has infinitely many branches in the complex plane, denoted by $W_k(z)$ for each integer $k\in\mathbb{Z}$. Corless et al. \cite[Eqn. (4.19)]{MR1414285} showed that, as $|z|\to\infty$,
\begin{equation}\label{WAE}
    W_k(z)=L_1-L_2+\frac{L_2}{L_1}+\frac{L_2(-2+L_2)}{2L_1^2}+O\bigg(\frac{L^3_2}{L^3_1}\bigg),
\end{equation}
where $L_1:=\log z+2\pi i k$ and $L_2:=\log L_1$, and $\log z$ is the principal logarithm.

Our starting point is the observation that $z_0$ can be described explicitly in terms of the Lambert $W$ function. We make this precise in Lemma 1. In what follows, let
\begin{equation}\label{lambda}
    \lambda(n):=\bigg(n+\frac{3}{2}\bigg)\frac{1}{z_0}.
\end{equation}
We now prove the following.
\begin{lemma}\label{L1}
    For all sufficiently large integers $n$,
\begin{equation}
    \begin{aligned}
        G(n)=-\bigg(n+\frac{3}{2}\bigg)&\bigg(\log \bigg(n+\frac{3}{2}\bigg)-\Re\log\lambda(n)-1+\Re\bigg(\frac{1}{\lambda(n)}\bigg)\bigg)\\&+\frac{1}{2}\log 2\pi+\frac{\Re\,\lambda(n)}{6n+9}+O\bigg(\frac{\log^3n}{n^3}\bigg)
    \end{aligned}
\end{equation}
and
\begin{equation}
    \begin{aligned}\label{F(n)}
        F(n)=\bigg(n+\frac{3}{2}\bigg)&\bigg(\Im\log \lambda(n)-\Im\bigg(\frac{1}{\lambda(n)}\bigg)\bigg)+\frac{\Im\lambda(n)}{6n+9}+O\bigg(\frac{\log^3n}{n^3}\bigg),
    \end{aligned}
\end{equation}
where $\lambda(n)$ is defined in \eqref{lambda} and satisfies
    \begin{equation}
        \lambda(n)=W_0(e^{T(n)})-\frac{W_0(e^{T(n)})^3}{\big(2\sqrt{3}n+3\sqrt{3}\big)^2\big(W_0(e^{T(n)})+1\big)}+O\bigg(\frac{\log^4 n}{n^4}\bigg).
    \end{equation}
\end{lemma}
\begin{proof}
    We begin with the saddle-point equation. Differentiating $\phi(z)$, we obtain
\begin{equation}\label{Derphi0}
    \phi'(z)=-\frac{n+1}{z}-\log (2\pi i)+\psi(z),
\end{equation}
where $\psi(z)$ is the digamma function. By Stirling's formula
\begin{equation}\label{Derphi}
    \phi'(z)=-\frac{n+1}{z}-\log (2\pi i)+\log z -\frac{1}{2z}-\frac{1}{12z^2}+O\bigg(\frac{1}{|z|^4}\bigg).
\end{equation}
Since $\phi'(z_0)=0$, setting $z=z_0$ in \eqref{Derphi} and multiplying by $z_0$ gives
\begin{equation}\label{SUB}
    z_0\log z_0-z_0\log(2\pi i)=n+\frac{3}{2}+\frac{1}{12z_0}+O\bigg(\frac{\log^3 n}{n^3}\bigg).
\end{equation}
We now use this relation to simplify $\phi(z_0)$. By Stirling's formula \cite[Eqn. (1.4.25)]{MR174795} and \eqref{SUB}, one has
\begin{equation}
   \begin{aligned}
        \phi(z_0)&=-(n+1)\log z_0-z_0\log(2\pi i)+\log\Gamma(z_0)\\&=-(n+1)\log z_0-z_0\log(2\pi  i)+\bigg(z_0-\frac{1}{2}\bigg)\log z_0-z_0+\frac{1}{2}\log 2 \pi+\frac{1}{12z_0}+O\bigg(\frac{\log^3 n}{n^3}\bigg)\\&=-\bigg(n+\frac{3}{2}\bigg)\log z_0+z_0\log z_0-z_0\log (2\pi i)-z_0+\frac{1}{2}\log 2\pi+\frac{1}{12z_0}+O\bigg(\frac{\log^3 n}{n^3}\bigg)\\&=-\bigg(n+\frac{3}{2}\bigg)\log z_0+\bigg(n+\frac{3}{2}\bigg)-z_0+\frac{1}{2}\log 2\pi+\frac{1}{6z_0}+O\bigg(\frac{\log^3 n}{n^3}\bigg).
   \end{aligned}
\end{equation}
As $z_0\in\Omega_n$, $\arg\lambda(n)=-\arg z_0\in(-\pi/4,0)$, and the principal logarithms satisfy
\begin{equation}\label{logz_0}
    \log z_0=\log\bigg(n+\frac{3}{2}\bigg)-\log\lambda(n).
\end{equation}
Therefore
\begin{equation}
    \begin{aligned}\label{REIM}
        \phi(z_0)=-\bigg(n+\frac{3}{2}\bigg)&\bigg(\log \bigg(n+\frac{3}{2}\bigg)-\log \lambda(n)-1+\frac{1}{\lambda(n)}\bigg)\\&+\frac{1}{2}\log 2\pi+\frac{\lambda(n)}{6n+9}+O\bigg(\frac{\log^3n}{n^3}\bigg).
    \end{aligned}
\end{equation}
Taking real and imaginary parts of \eqref{REIM} gives $G(n)$ and $F(n)$. To make this useful, however, we need an explicit description of $\lambda(n)$, which we now obtain.

Substituting \eqref{logz_0} in \eqref{SUB} gives
\begin{equation}\label{GE}
    \lambda(n)+\log \lambda(n)=T(n)-\frac{\lambda^2(n)}{\big(2\sqrt{3}n+3\sqrt{3}\big)^2}+O\bigg(\frac{\log^4n}{n^4}\bigg).
\end{equation}
Now consider the following equation obtained by dropping the last two terms on the right-hand side of \eqref{GE},
\begin{equation}\label{GE2}
    \tilde{\lambda}(n)+\log \tilde{\lambda}(n)=T(n).
\end{equation}
Unlike \eqref{GE}, this equation can be solved exactly. Exponentiating \eqref{GE2}, we obtain
\begin{equation}
    \tilde{\lambda}(n)e^{\tilde{\lambda}(n)}=e^{T(n)}.
\end{equation}
Hence every solution of \eqref{GE2} is of the form $W_k(e^{T(n)})$ for some $k\in\mathbb{Z}$. Since $\Im\,T(n)=-\pi/2$, it follows from the asymptotic expansion \eqref{WAE} that, for all sufficiently large integers $n$,
\begin{equation}
    W_0(e^{T(n)})+\log W_0(e^{T(n)})=T(n),
\end{equation}
so that $W_0(e^{T(n)})$ is a solution of \eqref{GE2}. Furthermore, it can also be checked that $\tilde{\lambda}(n):=W_0(e^{T(n)})$ is the unique branch of the Lambert $W$ function that solves \eqref{GE2} and satisfies $\tilde{\lambda}(n)=T(n)+O(\log\log n)$.

We now show that $\lambda(n)$ and $\tilde{\lambda}(n)$ are very close. We begin with a crude estimate. It is clear from the definition that $|\lambda(n)|\asymp \log n$ (cf. \eqref{a(n)}). Substituting this estimate in \eqref{GE} gives $\lambda(n)=T(n)+O(\log\log n)$. Therefore
\begin{equation}
    \log \lambda(n)=\log T(n)+O\bigg(\frac{\log\log n}{\log n}\bigg).
\end{equation}
Substituting this into \eqref{GE}, we obtain
\begin{equation}
    \lambda(n)=T(n)-\log T(n)+O\bigg(\frac{\log\log n}{\log n}\bigg).
\end{equation}
On the other hand, applying \eqref{WAE} with $k=0$ and $z=e^{T(n)}$, we see that $\tilde{\lambda}(n)$ admits the same expansion. Hence
\begin{equation}
\lambda(n)-\tilde{\lambda}(n)=O\bigg(\frac{\log\log n}{\log n}\bigg).
\end{equation}
Now, we shall use this crude estimate to obtain a sharper estimate for $\lambda(n)-\tilde{\lambda}(n)$. To this end, we define
\begin{equation}
   \delta(n):=\lambda(n)-\tilde{\lambda}(n)\quad\mathrm{and}\quad u(n):=\frac{\delta(n)}{\tilde{\lambda}(n)}.
\end{equation}
Since $|\tilde{\lambda}(n)|\asymp \log n$, the preceding estimate shows that $u(n)\to 0$ as $n\to\infty$. Moreover, $\lambda(n)=\tilde{\lambda}(n)(1+u(n))$, so that, for all sufficiently large $n$, $\log\lambda(n)=\log\tilde{\lambda}(n)+\log(1+u(n))$. Subtracting \eqref{GE2} from \eqref{GE} and using this identity, we obtain
\begin{equation}\label{GE3}
    \delta(n)+\log\big(1+u(n)\big)=-\frac{\lambda^2(n)}{\big(2\sqrt{3}n+3\sqrt{3}\big)^2}+O\bigg(\frac{\log^4n}{n^4}\bigg).
\end{equation}
We now estimate the left-hand side of the above equation. Since $u(n)\to0$ as $n\to\infty$, we may expand $\log\big(1+u(n)\big)$ in its power series to obtain
\begin{equation}\label{deltau}
    \begin{aligned}
         \delta(n)+\log\big(1+u(n)\big)=\delta(n)\bigg(1+\frac{1}{\tilde{\lambda}(n)}+O\bigg(\frac{|\delta(n)|}{|\tilde{\lambda}(n)|^2}\bigg)\bigg)=\delta(n)\bigg(1+O\bigg(\frac{1}{\log n}\bigg)\bigg).
    \end{aligned}
\end{equation}
Substituting the above estimate in \eqref{GE3} and rearranging the terms gives
\begin{equation}\label{delta}
    \delta(n)=O\bigg(\frac{\log^2 n}{n^2}\bigg).
\end{equation}
This suffices to prove that $\lambda(n)$ and $\tilde{\lambda}(n)$ are very close. 

Now we determine the leading term of $\delta(n)$. Since $|\lambda(n)|\asymp \log n$ and $|\tilde{\lambda}(n)|\asymp \log n$, it follows from \eqref{delta} that
\begin{equation}
    \lambda^2(n)-\tilde{\lambda}^2(n)=\big(\lambda(n)+\tilde{\lambda}(n)\big)\delta(n)=O\bigg(\frac{\log^3 n}{n^2}\bigg).
\end{equation}
Consequently,
\begin{equation}
    -\frac{\lambda^2(n)}{\big(2\sqrt{3}n+3\sqrt{3}\big)^2}=-\frac{\tilde{\lambda}^2(n)}{\big(2\sqrt{3}n+3\sqrt{3}\big)^2}+O\bigg(\frac{\log^3 n}{n^4}\bigg).
\end{equation}
Substituting this into \eqref{GE3} and using \eqref{deltau}, we obtain
\begin{equation}
    \delta(n)\bigg(1+\frac{1}{\tilde{\lambda}(n)}+O\bigg(\frac{|\delta(n)|}{|\tilde{\lambda}(n)|^2}\bigg)\bigg)=-\frac{\tilde{\lambda}^2(n)}{\big(2\sqrt{3}n+3\sqrt{3}\big)^2}+O\bigg(\frac{\log^4 n}{n^4}\bigg).
\end{equation}
Therefore,
\begin{equation}
    \lambda(n)=\tilde{\lambda}(n)-\frac{\tilde{\lambda}^3(n)}{\big(2\sqrt{3}n+3\sqrt{3}\big)^2\big(\tilde{\lambda}(n)+1\big)}+O\bigg(\frac{\log^4 n}{n^4}\bigg).
\end{equation}

The expressions for $\Re\,\phi(z_0)$ and $\Im\,\phi(z_0)$ can be calculated straightforwardly from \eqref{REIM}. This completes the proof of the lemma. 
\end{proof}

\subsection{Proof of Theorem 1}
Here we prove \autoref{Thm1} using \hyperref[L1]{Lemma \ref{L1}}. For convenience, we set $n':=n+3/2$ and
\begin{equation}\label{R}
    R(n):=\log\Bigg(\frac{n!e^{g(b)}}{\big(g''(b)^2+f''(b)^2\big)^{1/4}}\Bigg)=\log\bigg(\frac{n!e^{G(n)}}{|\phi''(z_0)|^{1/2}}\bigg).
\end{equation}
\begin{proof}
    Setting $N=0$ in \autoref{Matsuoka} and noting that $h_0=1$, we obtain
\begin{equation}
    \begin{aligned}\label{N=0}
        \gamma_n&=\sqrt{\frac{2}{\pi}}n!e^{g(b)}\big(g''(b)^2+f''(b)^2\big)^{-1/4}\cos\bigg\{f(b)-\frac{1}{2}\arctan\bigg(\frac{f''(b)}{g''(b)}\bigg)\bigg\}\\&\qquad\qquad\qquad\qquad+O\bigg(n!e^{g(b)}n^{1/6}\log^{-1/3}n\bigg)\\&=\sqrt{\frac{2}{\pi}}n!e^{g(b)}\big(g''(b)^2+f''(b)^2\big)^{-1/4}\bigg(\cos\bigg\{f(b)-\frac{1}{2}\arctan\bigg(\frac{f''(b)}{g''(b)}\bigg)\bigg\}\\&\qquad\qquad\qquad+O\bigg(\big(g''(b)^2+f''(b)^2\big)^{1/4}n^{1/6}\log^{-1/3}n\bigg)\bigg).
    \end{aligned}
\end{equation}
Since
\begin{equation}
    \frac{d^2}{dy^2}\phi(a+iy)=-\phi''(a+iy),
\end{equation}
taking real and imaginary parts at $y=b$ gives $g''(b)=-\Re\,\phi''(z_0)$ and $f''(b)=-\Im\,\phi''(z_0)$. Hence
\begin{equation}\label{3sq}
    g''(b)^2+f''(b)^2=\big(\Re\,\phi''(z_0)\big)^2+\big(\Im\,\phi''(z_0)\big)^2=|\phi''(z_0)|^2.
\end{equation}
Moreover, since $f''(b)/g''(b)=\Im\,\phi''(z_0)/\Re\,\phi''(z_0)$ and $\Re\,\phi''(z_0)=-g''(b)>0$ (cf. \eqref{g''}), we obtain
\begin{equation}
    -\frac{1}{2}\arctan\bigg(\frac{f''(b)}{g''(b)}\bigg)=-\frac{1}{2}\arg\phi''(z_0).
\end{equation}
Thus, both the error term and the phase of the cosine in \eqref{N=0} are determined by $\phi''(z_0)$. Therefore, we need an explicit expression for $\phi''(z_0)$, which we now derive. Differentiating \eqref{Derphi0} and using the asymptotic expansion (obtained by differentiating Eqn. (1.4.25) from \cite{MR174795})
\begin{equation}
    \psi'(z)=\frac{1}{z}+\frac{1}{2z^2}+\frac{1}{6z^3}+O\bigg(\frac{1}{|z|^5}\bigg),
\end{equation}
we obtain
\begin{equation}
   \begin{aligned}\label{phi''}
        \phi''(z_0)&=\frac{n'}{z_0^2}+\frac{1}{z_0} +\frac{1}{6z_0^3}+O\bigg(\frac{1}{|z_0|^5}\bigg)\\&=\frac{1}{n'}\bigg(\lambda^2(n)+\lambda(n)+\frac{\lambda^3(n)}{6n'^2}\bigg)+O\bigg(\frac{\log^5 n}{n^5}\bigg)\\&=\frac{1}{n'}\bigg(\lambda^2(n)+\lambda(n)+\frac{\lambda^3(n)}{6n'^2}\bigg)\bigg(1+O\bigg(\frac{\log^3 n}{n^4}\bigg)\bigg).
   \end{aligned}
\end{equation}
It follows that 
\begin{equation}
    |\phi''(z_0)|^{1/2}\asymp\frac{\log n}{\sqrt{n}}.
\end{equation}
By $\eqref{3sq}$, the error term in \eqref{N=0} is therefore 
\begin{equation}
    O\bigg(\frac{\log^{2/3}n}{n^{1/3}}\bigg).
\end{equation}
Moreover, taking arguments in \eqref{phi''}, and noting that $n'>0$, we get
\begin{equation}
    -\frac{1}{2}\arg\phi''(z_0)=-\frac{1}{2}\arg\bigg(\lambda^2(n)+\lambda(n)+\frac{\lambda^3(n)}{6n'^2}\bigg)+O\bigg(\frac{\log^3 n}{n^4}\bigg).
\end{equation}
Therefore
\begin{equation}
    \begin{aligned}\label{cos}
        \cos\bigg\{&f(b)-\frac{1}{2}\arctan\bigg(\frac{f''(b)}{g''(b)}\bigg)\bigg\}\\&=\cos\bigg\{F(n)-\frac{1}{2}\arg\bigg(\lambda^2(n)+\lambda(n)+\frac{\lambda^3(n)}{6n'^2}\bigg)+O\bigg(\frac{\log^3 n}{n^4}\bigg)\bigg\}.
    \end{aligned}
\end{equation}
Substituting \eqref{cos} in \eqref{N=0} and using \eqref{R} we obtain
\begin{equation}
    \begin{aligned}
        \gamma_n=\sqrt{\frac{2}{\pi}}e^{R(n)}&\bigg(\cos\bigg\{J(n)+O\bigg(\frac{\log^3 n}{n^4}\bigg)\bigg\}+O\bigg(\frac{\log^{2/3 }n}{n^{1/3}}\bigg)\bigg)
    \end{aligned}
\end{equation}
where
\begin{equation}
    J(n):=F(n)-\frac{1}{2}\arg\bigg(\lambda^2(n)+\lambda(n)+\frac{\lambda^3(n)}{6n'^2}\bigg).\label{J}
\end{equation}
Substituting \eqref{F(n)} for $F(n)$ gives the formula in the form as stated in \autoref{Thm1}. It remains to evaluate $R(n)$. Taking the logarithm of the absolute value of \eqref{phi''} we obtain
\begin{equation}\label{2phi''}
    -\frac{1}{2}\log|\phi''(z_0)|=\frac{1}{2}\log n'-\frac{1}{2}\log\bigg|\lambda^2(n)+\lambda(n)+\frac{\lambda^3(n)}{6n'^2}\bigg|+O\bigg(\frac{\log^3 n}{n^4}\bigg).
\end{equation}
Expanding the logarithm appearing in \eqref{R}, and using Stirling's formula for $\log n!$, \hyperref[L1]{Lemma \ref{L1}} for $G(n)$ and \eqref{2phi''} completes the proof of the theorem.
\end{proof}

\begin{figure}
    \centering
\includegraphics[width=0.9\linewidth]{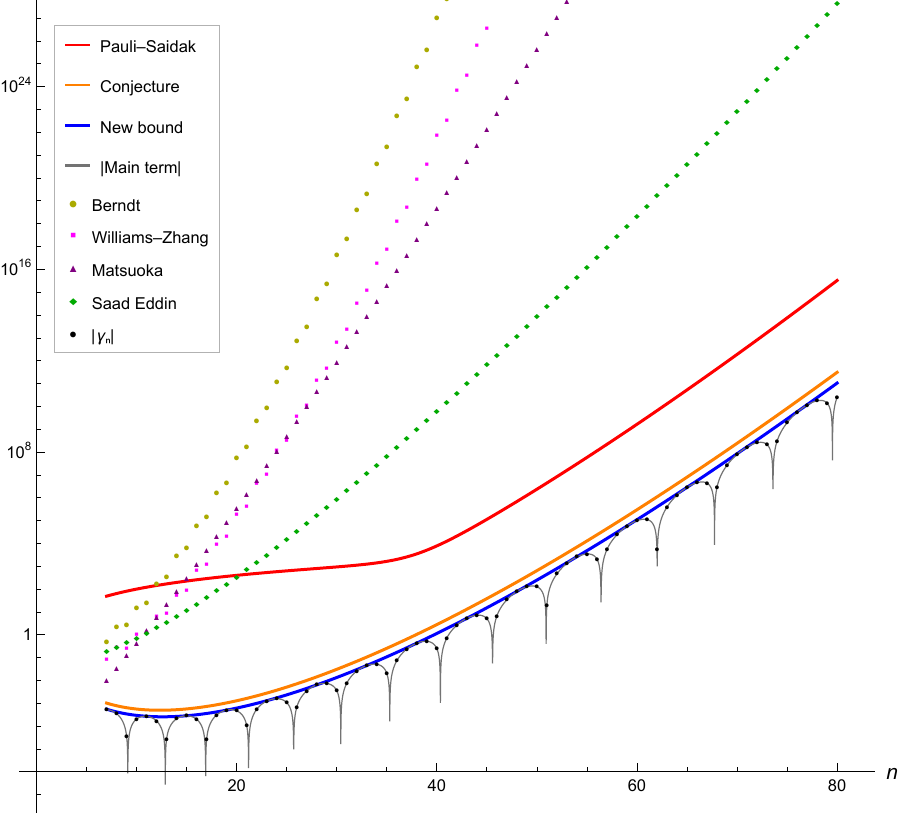}
    \caption{Comparison of the bounds for $|\gamma_n|$ on a logarithmic scale. The new bound (cf. \eqref{NB}) and the absolute value of the main term (cf. \eqref{MT}) are shown without the $O$-terms.  }
    \label{fig1}
\end{figure}

\begin{figure}
    \centering
\includegraphics[width=0.9\linewidth]{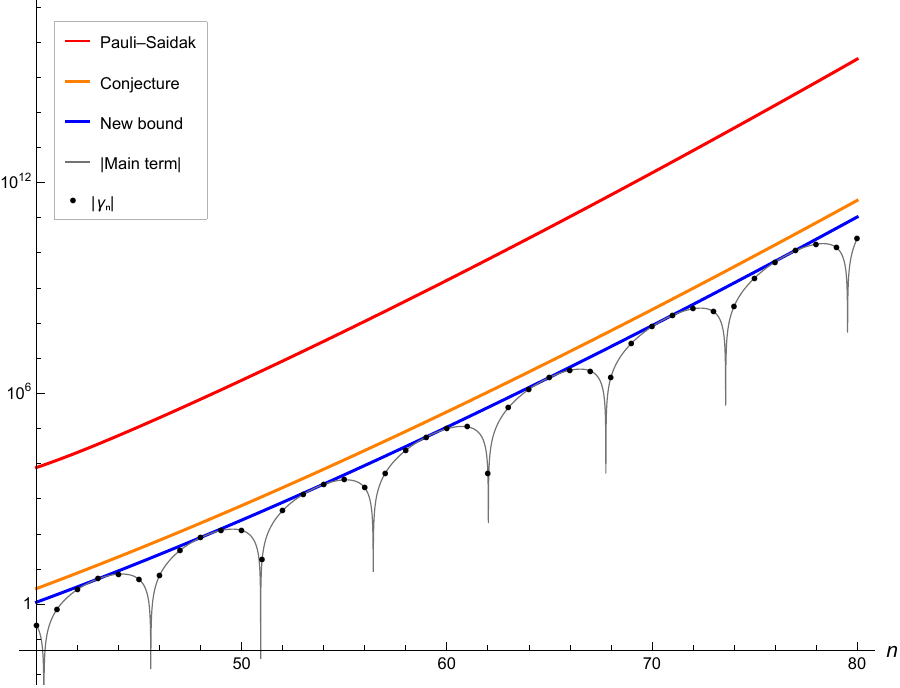}
    \caption{Part of \hyperref[fig1]{Figure \ref{fig1}}, enlarged.}
    \label{fig:placeholder}
\end{figure}

To put our results in context, we compare them with previously known bounds in \hyperref[fig1]{Figure \ref{fig1}} and \hyperref[tab:comparison]{Table \ref{tab:comparison}}. An enlarged part of \hyperref[fig1]{Figure \ref{fig1}} is shown in \hyperref[fig:placeholder]{Figure \ref{fig:placeholder}}.  The earliest of these bounds is due to Berndt \cite{MR289431}, who showed that
\begin{equation}
    |\gamma_n|\leq (3+(-1)^n)\frac{\Gamma(n)}{\pi^n}.
\end{equation}

\begin{table}[p]
\centering
\renewcommand{\arraystretch}{2.0}
\resizebox{\textwidth}{!}{%
\begin{tabular}{rlllll}
\hline
$n$ & Pauli--Saidak \cite{MR4671618} & Conjecture \cite{MR4280128} & \quad New bound & \quad Main term  & $\qquad\quad\gamma_n$  \\
\hline
10   & $1.0004832\times 10^{2}$    & $5.5954615\times 10^{-4}$   & $2.9608834\times 10^{-4}$   & $\phantom{-}2.0427384\times 10^{-4}$  & $\phantom{-}2.0533281\times 10^{-4}$ \\
100  & $2.7215231\times 10^{22}$   & $1.5759233\times 10^{18}$   & $4.9672125\times 10^{17}$   & $-4.2537006\times 10^{17}$            & $-4.2534016\times 10^{17}$ \\
150  & $4.5756730\times 10^{41}$   & $1.0823019\times 10^{37}$   & $2.9977502\times 10^{36}$   & $\phantom{-}8.0260268\times 10^{35}$  & $\phantom{-}8.0288537\times 10^{35}$ \\
200  & $6.5413861\times 10^{62}$   & $8.2307717\times 10^{57}$   & $2.0715004\times 10^{57}$   & $-6.9594225\times 10^{55}$            & $-6.9746497\times 10^{55}$ \\
250  & $2.2049090\times 10^{85}$   & $1.7038266\times 10^{80}$   & $3.9728549\times 10^{79}$   & $\phantom{-}3.0591544\times 10^{79}$  & $\phantom{-}3.0592129\times 10^{79}$ \\
300  & $8.4711748\times 10^{108}$  & $4.4005495\times 10^{103}$  & $9.6283526\times 10^{102}$  & $-5.5572672\times 10^{102}$           & $-5.5567282\times 10^{102}$ \\
350  & $2.3288119\times 10^{133}$  & $8.6541146\times 10^{127}$  & $1.7929004\times 10^{127}$  & $-1.7855058\times 10^{127}$           & $-1.7854461\times 10^{127}$ \\
400  & $3.3185429\times 10^{158}$  & $9.2313069\times 10^{152}$  & $1.8231091\times 10^{152}$  & $-1.7616774\times 10^{152}$           & $-1.7616422\times 10^{152}$ \\
450  & $1.9364490\times 10^{184}$  & $4.1740922\times 10^{178}$  & $7.8994930\times 10^{177}$  & $-7.7708787\times 10^{177}$           & $-7.7707021\times 10^{177}$ \\
500  & $3.8670886\times 10^{210}$  & $6.6374169\times 10^{204}$  & $1.2087589\times 10^{204}$  & $-1.1655454\times 10^{204}$           & $-1.1655053\times 10^{204}$ \\
550  & $2.2956228\times 10^{237}$  & $3.2071565\times 10^{231}$  & $5.6396148\times 10^{230}$  & $-2.8601579\times 10^{230}$           & $-2.8599679\times 10^{230}$ \\
600  & $3.6167943\times 10^{264}$  & $4.1881137\times 10^{258}$  & $7.1313982\times 10^{257}$  & $\phantom{-}3.5627037\times 10^{257}$ & $\phantom{-}3.5627463\times 10^{257}$ \\
650  & $1.3779573\times 10^{292}$  & $1.3427827\times 10^{286}$  & $2.2194091\times 10^{285}$  & $\phantom{-}2.1441173\times 10^{285}$ & $\phantom{-}2.1440540\times 10^{285}$ \\
700  & $1.1746180\times 10^{320}$  & $9.7579002\times 10^{313}$  & $1.5687831\times 10^{313}$  & $-3.5492612\times 10^{312}$           & $-3.5494522\times 10^{312}$ \\
750  & $2.0976985\times 10^{348}$  & $1.5021956\times 10^{342}$  & $2.3533596\times 10^{341}$  & $-2.1646512\times 10^{341}$           & $-2.1645871\times 10^{341}$ \\
800  & $7.4182482\times 10^{376}$  & $4.6239586\times 10^{370}$  & $7.0699232\times 10^{369}$  & $\phantom{-}4.9135696\times 10^{369}$ & $\phantom{-}4.9135406\times 10^{369}$ \\
850  & $4.9474668\times 10^{405}$  & $2.7071767\times 10^{399}$  & $4.0453884\times 10^{398}$  & $\phantom{-}1.0481917\times 10^{398}$ & $\phantom{-}1.0481137\times 10^{398}$ \\
900  & $5.9631473\times 10^{434}$  & $2.8860073\times 10^{428}$  & $4.2200871\times 10^{427}$  & $-3.8179450\times 10^{427}$           & $-3.8178431\times 10^{427}$ \\
950  & $1.2510528\times 10^{464}$  & $5.3913444\times 10^{457}$  & $7.7229453\times 10^{456}$  & $\phantom{-}7.4263377\times 10^{456}$ & $\phantom{-}7.4262177\times 10^{456}$ \\
1000 & $4.4189788\times 10^{493}$  & $1.7059006\times 10^{487}$  & $2.3962680\times 10^{486}$  & $-1.5709625\times 10^{486}$           & $-1.5709538\times 10^{486}$ \\
1200 & $7.2525208\times 10^{613}$  & $1.8942805\times 10^{607}$  & $2.4808988\times 10^{606}$  & $\phantom{-}8.6844007\times 10^{605}$ & $\phantom{-}8.6840299\times 10^{605}$ \\
1300 & $1.1080908\times 10^{675}$  & $2.4385461\times 10^{668}$  & $3.0963636\times 10^{667}$  & $\phantom{-}7.4827575\times 10^{666}$ & $\phantom{-}7.4823742\times 10^{666}$ \\
1400 & $7.2676255\times 10^{736}$  & $1.3649386\times 10^{730}$  & $1.6839905\times 10^{729}$  & $-4.0972199\times 10^{728}$           & $-4.0972873\times 10^{728}$ \\
1500 & $1.8177550\times 10^{799}$  & $2.9458655\times 10^{792}$  & $3.5380770\times 10^{791}$  & $-3.1170517\times 10^{791}$           & $-3.1170181\times 10^{791}$ \\
1600 & $1.5665569\times 10^{862}$  & $2.2118235\times 10^{855}$  & $2.5903112\times 10^{854}$  & $-2.0574934\times 10^{854}$           & $-2.0574536\times 10^{854}$ \\
1700 & $4.2605893\times 10^{925}$  & $5.2852070\times 10^{918}$  & $6.0442885\times 10^{917}$  & $\phantom{-}2.9861798\times 10^{917}$ & $\phantom{-}2.9861730\times 10^{917}$ \\
1800 & $3.3865741\times 10^{989}$  & $3.7186154\times 10^{982}$  & $4.1582467\times 10^{981}$  & $\phantom{-}3.3153704\times 10^{981}$ & $\phantom{-}3.3153125\times 10^{981}$ \\
1900 & $7.3526318\times 10^{1053}$ & $7.1941565\times 10^{1046}$ & $7.8751616\times 10^{1045}$ & $-6.7431291\times 10^{1045}$          & $-6.7430678\times 10^{1045}$ \\
2000 & $4.1063168\times 10^{1118}$ & $3.6016046\times 10^{1111}$ & $3.8634980\times 10^{1110}$ & $\phantom{-}2.6802539\times 10^{1109}$ & $\phantom{-}2.6804247\times 10^{1109}$ \\
\hline
\end{tabular}}
\caption{Comparison of our results with the bound of Pauli and Saidak, the conjectured bound and the true values of $\gamma_n$. The fourth and fifth columns are the main terms of \eqref{NB} and \eqref{MT}, evaluated without the $O$-terms, whose implied constants are not explicit.}
\label{tab:comparison}
\end{table}

Williams and Zhang \cite{MR1302533} later obtained 
\begin{equation}
    |\gamma_n|\leq(3+(-1)^n)\frac{\Gamma(2n+1)}{n^{n+1}(2\pi)^n}.
\end{equation}
Matsuoka's bound \cite{MR827790}, already mentioned in the introduction, states that $|\gamma_n|<10^{-4}(\log n)^n$. A sharper bound was given by Saad Eddin \cite{MR2997785}: if we put 
\begin{equation}
    \theta_n:=\frac{n+1}{\log\big(\frac{2(n+1)}{\pi}\big)}-1,
\end{equation}
then
\begin{equation}
    |\gamma_n|\leq 2\sqrt{2}\Gamma(n+1)e^{-(n+1)\log\theta_n+\theta_n\log\big(\frac{2\theta_n}{\pi e}\big)}\bigg(1+2^{-\theta_n-1}\frac{\theta_n+1}{\theta_n-1}\bigg).
\end{equation}
Finally, we include the bounds \eqref{PS} and \eqref{CON}, both of which are expressed in terms of $W_0\big(i\frac{n}{2\pi}\big)$.

\subsection*{Declaration and Acknowledgment}

 The authors declare that no funding was received during the preparation of this manuscript. The authors declare no competing interests. No data were generated or analyzed in this study.

The authors used AI tools to assist with literature searches, grammar checking, bibliography management, LaTeX assistance and proofreading. All the main mathematical ideas are the authors' own, and the authors take full responsible for the content of the paper.

\clearpage
\bibliographystyle{amsplain}
\bibliography{references}
\end{document}